\documentclass{amsart}

\usepackage{amsfonts, amsmath, amssymb, amsthm, amsxtra, latexsym}

\theoremstyle{plain}

\newtheorem{theorem}{Theorem}[section]

\theoremstyle{definition}

\newtheorem{remark}[theorem]{Remark}

\newcommand{\R}{{\mathbb R}}
\newcommand{\C}{{\mathbb C}}

\numberwithin{equation}{section}

\begin{document}


\title[Remarks on a proof of Paley--Wiener theorems]{Some remarks on a proof of geometrical Paley--Wiener theorems for the Dunkl transform}


\translator{}

\dedicatory{}

\author{Marcel de Jeu}

\begin{abstract}
We argue that a proof of the geometrical form of the Paley--Wiener theorems for
the Dunkl transform in the literature is not correct.
\end{abstract}

\date{}



\address{M.F.E.~de~Jeu\\
         Mathematical Institute\\
         Leiden University\\
         P.O. Box 9512\\
         2300 RA Leiden\\
         The Netherlands}



\email{mdejeu@math.leidenuniv.nl}

\urladdr{}

\maketitle


\section{Introduction}

In \cite{trimeche} a proof of the geometrical form of the Paley--Wiener
theorems for the Dunkl transform is presented. In our opinion, however, this
proof is not correct. We have informed the author of the details underlying
this opinion from November 2003 onward, but at the time of writing he has not
refuted our remarks, or agreed with them, or given an alternative correct
proof, or expressed the intention to publish an Erratum.

The material which was communicated to the author is presented below. It is our
opinion that at this moment the geometrical forms of the Paley--Wiener theorems
for the Dunkl transform are still unproven.

\section{Arguments}

In \cite{trimeche} the geometric form of the Paley-Wiener theorem for the Dunkl
transform is stated for functions and for distributions, as Theorems~6.2 and
6.3, respectively. The crucial ingredient in the proof of these results is
Proposition~6.3. Our arguments concern the proof of this proposition. In
formulating them, we will use the notation and definitions of \cite{trimeche}.

In the proof of Proposition~6.3, a $W$-invariant compact convex subset $E$ of
$\R^d$ is considered. For $x\in E$ fixed, the function $f_x$ on $\R^d$ is defined as
\begin{equation}\label{eq:fxdefinition}
f_x(y)=\frac{e^{-i(x,y)}}{(1+\Vert y \Vert^2)^p}\quad(y\in\R^d),
\end{equation}
where $p$ is an integer such that $p\geq\gamma+d/2+1$. Since the constant
$\gamma$ is assumed to be strictly positive in line 3 on page 29, we see that
$p\geq 2$.

The function $F_x$ is defined on $\R^d$ in equation (69) as essentially the
inverse Dunkl transform of $f_x$, namely
\[
F_x(t)=\int_{\R^d} f_x(y)K(iy,t)\omega_k(y)dy\quad(t\in\R^d).
\]
Following this definition, it is observed that $F_x$ is continuous.

After a computation involving Riemann sums and contour integration, it is then
concluded on line -3 on page 31 that $F_x$ has support in the set $E$. Since
$E$ is compact, one sees---if the line of reasoning in \cite{trimeche} is
correct---that $F_x$ is a compactly supported continuous function.

But this never holds. Indeed, if $F_x$ were compactly supported, then $f_x$
could be reconstructed from $F_x$ as is stated in the last line on page 31 (an
application of the inversion theorem for the Dunkl transform):
\begin{equation}\label{eq:reconstruction}
f_x(y)=\frac{C_k^2}{2^{2\gamma+d}}\int_{\R^d}
F_x(t)K(-iy,t)\omega_k(t)\,dt\quad(y\in\R^d).
\end{equation}
However, as with the ordinary Fourier transform, the fact that $F_x$ is a
continuous function with compact support implies \cite[Part 3 of Lemma
4.4]{dunkltransform} that its Dunkl transform has an entire extension to
$\C^d$, i.e., that $f_x$ has an entire extension to $\C^d$. But it follows from
\eqref{eq:fxdefinition}, written as
\[
f_x(y)=\frac{e^{-i(x,y)}}{(1+(y,y))^p},
\]
where $(\,.\,,\,.\,)$ is the holomorphic standard bilinear form on $\C^d$, that
$f_x$ has no such extension from $\R^d$ to $\C^d$, since this extension would
have a pole along the divisor $\{y\in\C^d\mid (y,y)=-1\}$, as a consequence of
the fact that $p\geq 2$. This is a contradiction, so $F_x$ can never have
compact support.

The above argument shows that $F_x$ can never be a continuous function with
compact support, but apart from that it is also easy to give an elementary
counterexample to the claim in \cite{trimeche} that $F_x$ always has these
properties, as follows. Consider the $W$-invariant compact convex set $E=\{0\}$
and choose $x=0\in E$. Then the fact that the continuous function $F_0$ has
support in $E$ implies that $F_0=0$. But in that case \eqref{eq:reconstruction}
implies that $f_0=0$, contradicting \eqref{eq:fxdefinition}.

\medskip
The fact that $F_x$ has support in $E$ is the cornerstone of the proof of
Proposition~6.3 and therefore also of the proof of the geometrical
Paley--Wiener theorems in \cite{trimeche}. Since this statement is in fact
false, as shown above, we consider the proof of the geometrical Paley--Wiener
theorems in \cite{trimeche} to be incorrect.

\begin{remark}
The crucial mistake in the line of reasoning which leads to the statement that
$F_x$ has support in $E$ appears to be the following.

In line -9 on page 31 it is stated that, after shifting the domain of
integration over $i\eta$ for $\eta\in\R^d$, one has
\begin{equation}\label{eq:shiftedintegral}
F_x(t)=\frac{1}{(i\pi)^d}\int_{\R^d}
e^{-i(x,u+i\eta)}K(i(u+i\eta),t)m_k(u+i\eta) \,du \quad(t\in\R^d).
\end{equation}
Following this it is claimed that an application of equation (20) in
\cite{trimeche}, i.e., an application of the estimate
\begin{equation}\label{eq:estimates}
|K(x,z)|\leq \exp\left(\Vert x\Vert \, \Vert \textup{Re}\, z\Vert \right)
\quad(x\in\R^d,\,z\in \C^d)
\end{equation}
implies---here one also uses the symmetry of the Dunkl kernel---that
\begin{equation}\label{eq:incorrectestimate}
|F_x(t)|\leq\frac{e^{(I_E(\eta)-(t, \eta))}}{\pi^d}\int_{\R^d} |m_k(u+i\eta)|
\,du \quad(t\in\R^d),
\end{equation}
where $I_E(\eta)=\sup_{x\in E}(x,\eta)$. Next it is argued that
\eqref{eq:incorrectestimate} implies the crucial fact that $F_x(t)=0$ if
$t\notin E$, by shifting the domain to infinity in a direction depending on
$t$.

However, the estimate in \eqref{eq:incorrectestimate} does not follow from
\eqref{eq:shiftedintegral} and \eqref{eq:estimates}.  One can only conclude
that
\begin{equation}\label{eq:correctestimate}
|F_x(t)|\leq\frac{e^{(I_E(\eta)+\Vert t\Vert
\Vert\eta\Vert)}}{\pi^d}\int_{\R^d} |m_k(u+i\eta)| \,du \quad(t\in\R^d).
\end{equation}
But then shifting the domain of integration no longer works. If ones chooses,
e.g., for $E$ a ball with radius $R$ centered at the origin, then a factor
$e^{(R+\Vert t\Vert)\Vert\eta\Vert}$ appears before the integral in
\eqref{eq:correctestimate}, and this factor has exponential growth in any
direction for $\eta$, rather than exponential decay.
\end{remark}

\begin{remark}
If the multiplicity is strictly positive, then in the one-dimensional case the
kernel $K(iy,x)$ has exponential growth if $0\neq x\in\R$ is fixed and $y$
tends to infinity in the imaginary direction through the upper or lower half
plane. The choice of the half plane is immaterial, quite in contrast to the
ordinary exponential function. The available estimates for the Dunkl kernel
suggest that the same behaviour may occur in any dimension. As a consequence,
shifting the domain of integration is in our opinion not likely to work as long
as the integrand contains the Dunkl kernel as a factor. This technical
difficulty was already observed in \cite{thesis} and the above remarks about
the proof of the geometrical Paley--Wiener theorems in \cite{trimeche} seem to
underline this observation.
\end{remark}


\end{document}